\documentclass[a4paper]{amsart}
\usepackage{amssymb,latexsym,amsmath}
\usepackage{amsmath}

\theoremstyle{plain}
\newtheorem{theorem}{Theorem}
\newtheorem{lemma}{Lemma}
\newtheorem{corollary}{Corollary}

\theoremstyle{definition}

\newtheorem{example}{Example}

\begin{document}
\title{Normability of Probabilistic Normed Spaces}
\author[B. Lafuerza--Guill\'en]{Bernardo Lafuerza--Guill\'en}
\address{Departamento de Estad\'\i stica y Matem\'atica Aplicada\\
Universidad de Almer\'\i a\\
04120 Almer\'\i a, Spain} \email{blafuerz@ual.es}
\author[J.A. Rodr\'\i guez--Lallena]{Jos\'e Antonio Rodr\'\i guez--Lallena}
\email{jarodrig@ual.es}
\author[C. Sempi]{Carlo Sempi}
\address{Dipartimento di Matematica \lq\lq Ennio De Giorgi\rq\rq\\
Universit\`a di Lecce\\
73100 Lecce, Italy} \email{carlo.sempi@unile.it}

\date{\today}
\begin{abstract} Relying on Kolmogorov's classical characterization
of norm\-able Topological Vector spaces, we study the normability
of those Probabilistic Normed Spaces that are also Topological
Vector spaces and provide a characterization of normable \v
Serstnev spaces. We also study the normability of other two
classes of Probabilistic Normed Spaces.
\end{abstract}

\maketitle

\newcommand{\al}{\alpha}
\newcommand{\be}{\beta}
\newcommand{\de}{\delta}
\newcommand{\De}{\Delta}
\newcommand{\ep}{\epsilon}
\newcommand{\la}{\lambda}
\newcommand{\tet}{\theta}
\newcommand{\D}{\mathcal D}
\newcommand{\F}{\mathcal F}
\newcommand{\R}{\mathbf R}
\newcommand{\N}{\mathbf N}
\newcommand{\lp}{\left(}
\newcommand{\rp}{\right)}
\newcommand{\lsp}{\left[}
\newcommand{\rsp}{\right]}
\newcommand{\lop}{\left]}
\newcommand{\rop}{\right[}
\newcommand{\xra}{\xrightarrow}
\newcommand{\sm}{\setminus}
\newcommand{\ie}{i.e.}

\section{Introduction}

Probabilistic Normed spaces were introduced by \v Serstnev in
\cite{anS63}; their definition was generalized in \cite{ASS93}, a
paper that revived the study of these spaces. We recall the
definition, the properties and the examples of Probabilistic
Normed spaces that will be used in the following.

Let $\De$ be the space of distribution functions and
$\De^+:=\{F\in\De:F(0)=0\}$ the subset of distance distribution
functions \cite{book}. The space $\De$ can be metrized in several
equivalent ways \cite{dS71,bS75,cS82,mT85} in such a manner that
the metric topology coincides with the topology of weak
convergence for distribution functions. Here, we assume that $\De$
is metrized by the Sibley metric $d_S$, which is the same metric
denoted by $d_L$ in \cite{book}. We shall also consider the subset
$\D^+\subset\De^+$ of the proper distance distribution functions,
\ie\ those $F\in\De^+$ for which $\lim_{x\to +\infty} F(x)=1$.

A \textit{triangle function} is a mapping
$\tau:\De^+\times\De^+\to\De^+$ that is commutative, associative,
nondecreasing in each variable and has $\ep_0$ as identity, where
$\ep_a$ $(a\le +\infty)$ is the distribution function defined by
\[
\ep_a(t):=\begin{cases} 0, &\qquad t\le a,\\
1, &\qquad t>a. \end{cases}
\]
Given a nonempty set $S$, a mapping $\F$ from $S\times S$ into
$\De^+$ and a triangle function $\tau$, a \textit{Probabilistic
Metric Space} (briefly a PM space) is the triple $(S,\F,\tau)$
with the following properties, where we set $F_{p,q}:=\F(p,q)$,
\begin{enumerate}
\item[(M1)] $F_{p,q}=\ep_0$ if, and only if, $p=q$; \item[(M2)]
$F_{p,q}=F_{q,p}$ for all $p$ and $q$ in $S$; \item[(M3)]
$F_{p,r}\ge\tau\lp F_{p,q},F_{q,r}\rp$ for all $p,\, q,\,r\in S$.
\end{enumerate}

A \textit{Probabilistic Normed Space} (briefly a PN space) is a
quadruple $(V,\nu,\tau,\tau^*)$, where $V$ is a vector space,
$\tau$ and $\tau^*$ are continuous triangle functions such that
$\tau\le\tau^*$ and $\nu$ is a mapping from $V$ into $\De^+$,
called the \textit{probabilistic norm}, such that for every choice
of $p$ and $q$ in $V$ the following conditions hold:
\begin{enumerate}
\item[(N1)] $\nu_p=\ep_0$ if, and only if, $p=\tet$ ($\tet$ is the
null vector in $V$); \item[(N2)] $\nu_{-p}=\nu_p$; \item[(N3)]
$\nu_{p+q}\ge\tau\lp\nu_p,\nu_q\rp$; \item[(N4)]
$\nu_p\le\tau^*\lp\nu_{\la p},\nu_{(1-\la)p}\rp$ for every
$\la\in\lsp 0,1\rsp$.
\end{enumerate}
The pair $(V,\nu)$ is called a \textit{Probabilistic Seminormed
space} (PSN space for short) if $\nu$ satifies (N1) and (N2).

There are special PN spaces, only some of which we list below; for
the others we refer to \cite{LRS97}.

When there is a continuous $t$--norm $T$ (see \cite{book, KMP})
such that $\tau=\tau_T$ and $\tau^*=\tau_{T^*}$, where
$T^*(x,y):=1-T(1-x,1-y)$,
\[
\tau_T(F,G)(x):=\sup_{s+t=x} T\lp F(s),G(t)\rp\ \ \mbox{and}\ \
\tau_{T^*}(F,G)(x):=\inf_{s+t=x} T^*\lp F(s),G(t)\rp
\]
the PN space $(V,\nu,\tau_T,\tau_{T^*})$ is called a
\textit{Menger} PN space, and is denoted by $(V,\nu,T)$.

A PN space is called a \textit{\v{S}erstnev space} if it satisfies
(N1), (N3) and the following condition, which implies both (N2)
and (N4)
\begin{equation*}
\forall\ p\in V\ \ \forall\ \al\in\R\sm\{0\}\ \ \forall\ x>0\qquad
\nu_{\al p}(x)=\nu_p\lp\frac{x}{|\al|}\rp.
\end{equation*}

One speaks of an \textit{equilateral} PN space when there is
$F\in\De^+$ different from both $\ep_0$ and $\ep_{\infty}$ such
that, for every $p\ne\tet$, $\nu_p=F$, and when
$\tau=\tau^*=\mathbf M$, which is the triangle function defined
for $G$ and $H$ in $\De^+$ by $\mathbf
M(G,H)(x):=\min\{G(x),H(x)\}$.

Let $G\in\De^+$ be different from $\ep_0$ and from $\ep_{\infty}$
and let $(V,\|\cdot\|)$ be a normed space; then define, for
$p\ne\tet$,
\[
\nu_p(x):=G\lp\frac{x}{\| p\,\|}\rp.
\]
Then $(V,\nu,M)$ is a Menger space denoted by $(V,\|\cdot\|,G,M)$
($M(x,y):=\min\{x,y\}$). In the same conditions, if $\nu$ is
defined by
\[
\nu_p(x):=G\lp\frac{x}{\| p\,\|^{\al}}\rp,
\]
with $\al\ge 0$, then the pair $(V,\nu)$ is a PSN space called
\textit{$\al$--simple} and it is denoted by $(V,\|\cdot\|,G;\al)$.
The $\al$-simple spaces can be endowed with a structure of PN
space in a very general setting (when $\al>1$, $G$ should be a
continuous and strictly increasing function in $\D^+$: see
\cite{LRS97}).

See \cite{LRS95,LRS97,LRS99} for properties of PN spaces.

If $(V,\nu,\tau,\tau^*)$ is a PN space, a mapping $\mathcal
F\colon V\times V\to \De^+$ can be defined through
\begin{equation}\label{E:pd}
\mathcal F(p,q):=\nu_{p-q}.
\end{equation}
This function $\mathcal F$ makes $(V,\mathcal F,\tau)$ a
Probabilistic Metric Space. Every PM space can be endowed with the
strong topology; this topology is generated by the \textit{strong
neighbourhoods}, which are defined as follows: for every $t>0$,
the neighbourhood $N_p(t)$ at a point $p$ of $V$ is defined by
\[
N_p(t):=\left\{q\in V:d_S(\nu_{p-q},\ep_0)<t\right\} =\left\{q\in
V:\nu_{p-q}(t)>1-t\right\}.
\]
It is known (see \cite{book}) that $(V,\mathcal F,\tau)$, where
$\mathcal F$ is defined by \eqref{E:pd}, and therefore
$(V,\nu,\tau,\tau^*)$, is a Hausdorff space, and hence, a $T_1$
space; moreover, it is metrizable.

\section{PN spaces and TV spaces}

A result from \cite{ASS97} can be rephrased for the purpose of the
present paper in the following form

\begin{theorem}[\textbf{Alsina, Schweizer, Sklar}]\label{T:cont} Every \textrm{PN} space
$(V,\nu,\tau,\tau^*)$, when it is endowed with the strong topology
induced by the probabilistic norm $\nu$, is a topological vector
space if, and only if, for every $p\in V$ the map from $\R$ into
$V$ defined by
\begin{equation}\label{E:cont}
\al\mapsto \al\,p
\end{equation}
is continuous.
\end{theorem}

It was proved in \cite[Theorem 4]{ASS97} that, if the triangle
function $\tau^*$ is Archimedean, \ie\ if $\tau^*$ admits no
idempotents other than $\ep_0$ and $\ep_{\infty}$ (\cite{book}),
then the mapping \eqref{E:cont} is continuous and, as a
consequence, the PN space $(V,\nu,\tau,\tau^*)$ is a Topological
Vector space (=TV space).

\begin{theorem}\label{T:special} $(a)$ No equilateral space $(V,F,\mathbf M)$ is a\/
\textrm{TV} space.

$(b)$ A \v Serstnev space $(V,\nu,\tau)$ is a\/ \textrm{TV} space
if, and only if, the probabilistic norm $\nu$ maps $V$ into $\D^+$
rather than into $\De^+$, viz. $\nu(V)\subseteq\D^+$.

$(c)$ If $G$ is a continuous and strictly increasing function in
$\D^+$, then the $\al$--simple space $(V,\|\cdot\|,G;\al)$ is a
\textrm{TV} space.
\end{theorem}
\begin{proof} Let $\tet$ denote the null vector of the linear space $V$.
Since any PM space and, hence, any PN space, can be metrized, one
can limit oneself to investigate the behaviour of sequences.
Moreover, because of the linear structure of $V$, one can take
$p\ne\tet$ and an arbitrary sequence $\{\al_n\}$ with $\al_n\ne 0$
($n\in\N$) such that $\al_n\to 0$ as $n$ tends to $+\infty$.

(a) For every $n\in\N$, one has $\nu_{\al_n p}=F$ while
$\nu_{\tet}=\ep_0$. Therefore the map (2) is not continuous.

(b) If $\nu$ maps $V$ into $\D^+$, then, for every $t>0$, one has
\[
\nu_{\al_n p}(t)=\nu_p\lp \frac t{|\al_n|}\rp \xrightarrow[n\to
+\infty] {}1,
\]
whence the assertion. Conversely, if there exists at least one
$p\in V$ such that $\nu_p\in\De^+\sm\mathcal{D}^+$, namely such
that $\nu_p(x)\xra[x\to +\infty]{}\be<1$, then
\[
\nu_{\al_n p}(x)=\nu_p\lp\frac x{|\al_n|}\rp\xra[n\to
+\infty]{}\be<1,
\]
so that the mapping $\al\mapsto\al p$ is not continuous.

(c) Let $\{\al_n\}$ be a sequence of real numbers that goes to 0.
Then, for all $p\in V$ and $x>0$, one has
\[
\nu_{\al_n\,p}(x)=G\lp\frac x{\|\al_n\,p\,\|^{\al}}\rp =G\lp\frac
x{|\al_n|^{\al}\,\|p\,\|^{\al}}\rp.
\]
And $\lim_{n\to +\infty} \nu_{\al_n\,p}(x)=1$ since $G$ belongs to
$\D^+$.
\end{proof}

\begin{corollary}\label{C:cor1}
$(a)$ A simple space $(V,\|\cdot\|,G,M)$ is a\/ \textrm{TV} space
if, and only if, $G$ belongs to $\D^+$.

$(b)$ An\/ \textrm{EN }space $(S,\nu)$ is a\/ \textrm{TV} space
if, and only if, $\nu$ maps $S$ into $\D^+$, \ie\ $\nu_p\in\D^+$
for every $p\in S$.
\end{corollary}
\begin{proof} Both simple spaces and EN spaces are \v Serstnev spaces
(see \cite{LRS97}).
\end{proof}

If $\nu_p(x)$ is viewed as the probability $P(\|p\,\|<x)$ that the
usual norm of $p$ is less than $x$, then, the fact that, for some
$p\in V$, $\nu_p$ does not belong to $\D^+$ means that
$P(\|p\,\|<+\infty)<1$; this is to be regarded as being \lq\lq
odd\rq\rq. Therefore we shall call \textit{characteristic} any PN
space $(V,\nu,\tau,\tau^*)$ such that $\nu(V)\subseteq\D^+$, or,
equivalently, such that $\nu_p$ belongs to $\D^+$ for every $p\in
V$. Thus Theorem \ref{T:special} (b) and (c) can be rephrased as
follows.

\begin{theorem}\label{T:TVSS} $(a)$ A \v Serstnev space $(V,\nu,\tau)$ is
a\/ \textrm{TV} space if, and only if, it is characteristic.

$(b)$ Let $G$ be a continuous and strictly increasing distribution
function in $\De^+$. Then, the $\al$--simple space
$(V,\|\cdot\|,G;\al)$ is a\/ \textrm{TV} space if, and only if, it
is characteristic.
\end{theorem}

However, in general PN spaces, the condition $\nu(V)\subseteq\D^+$
is not necessary to obtain a TV space: see Theorem \ref{T:f}
below.

\section{Normability of PN spaces}

If $(V,\nu,\tau,\tau^*)$ is a TV space, the question naturally
arises of whether it is also normable; in other words, whether
there is a norm on $V$ that generates the strong topology. This
question had been broached by Prochaska (\cite{bP67}) in the case
of \v Serstnev PN spaces. For this case, we shall provide a
complete characterization of those characteristic \v Serstnev PN
spaces that are indeed normable (see Theorem \ref{T:ch*S} further
on). In the process, we shall need Kolmogorov's classical
characterization of normability for $T_1$ spaces (\cite{anK34}).

\begin{theorem}[\textbf{Kolmogorov}]\label{T:kol} A $T_1$\/ \textrm{TV} space is normable
if, and only if, there is a neighbourhood of the origin $\tet$
that is convex and topologically bounded.
\end{theorem}

Here, we have called \textit{topologically bounded\/} a set $A$ in
a TV space $E$ when, for every sequence $\{\al_n\}$ of real
numbers that converges to $0$ as $n$ tends to $+\infty$ and for
every sequence $\{p_n\}$ of elements of $A$, one has
$\al_n\,p_n\to \tet$ in the topology of $E$.

We recall that the \textit{probabilistic radius} of a set $A$ in a
PN space $(V,\nu,\tau,\tau^*)$ is the distance distribution
function $R_A$ given by $R_A(x):=\Phi_A(x-)$ ($=\lim_{u\to
x^-}\Phi_A(u)$) for all $x \in ]0,\infty[$, where
$\Phi_A(u):=\inf\{\nu_p(u):p\in A\}$ for all $u \in ]0,\infty[$
(see \cite{LRS99}). A subset $A$ of a PN space
$(V,\nu,\tau,\tau^*)$ is said to be $\D$--\textit{bounded} if, and
only if, there exists a distribution function $G\in\D^+$ such that
$\nu_p\ge G$ for every $p\in A$: in fact, we can take $G=R_A$.

\subsection{The case of \v Serstnev spaces}

In characterizing normable \v Serstnev spaces we shall need the
following result.

\begin{theorem}\label{T:eq*S} In a characteristic \v Serstnev space $(V,\nu,\tau)$
the following statements are equivalent for a subset $A$ of $V$:
\begin{enumerate}
\item[(a)] $A$ is $\D$--bounded; \item[(b)] $A$ is topologically
bounded.
\end{enumerate}
\end{theorem}
\begin{proof} (a) $\Longrightarrow$ (b) Let $A$ any $\D$--bounded subset
of $V$ and let $\{p_n\}$ be any sequence of elements of $A$ and
$\{\al_n\}$ any sequence of real numbers that converges to $0$;
there is no loss of generality in assuming $\al_n\ne 0$ for every
$n\in\N$. Then, for every $x>0$, and for every $n\in\N$,
\[
\nu_{\al_n p_n}(x)=\nu_{p_n}\lp \frac x{|\al_n|}\rp \ge R_A\lp
\frac x{|\al_n|}\rp \xrightarrow[n\to +\infty]{}1.
\]
Thus $\al_n\,p_n\to\tet$ in the strong topology and $A$ is
topologically bounded.

(b) $\Longrightarrow$ (a) Let $A$ be a subset of $V$ which is not
$\D$--bounded. Then
\[
R_A(x) \xrightarrow[x\to +\infty]{} \gamma <1.
\]
By definition of $R_A$, for every $n\in\N$ there is $p_n\in A$
such that
\[
\nu_{p_n}(n^2)<\frac {1+\gamma}2<1.
\]
If $\al_n=1/n$, then
\[
\nu_{\al_n\,p_n}(1/2)\le \nu_{\al_n\,p_n}(n)=\nu_{p_n}(n^2) <\frac
{1+\gamma}2 <1,
\]
which shows that $\{\nu_{\al_n\,p_n}\}$ does not tend to $\ep_0$,
even if it has a weak limit, viz. $\{\al_n\,p_n\}$ does not tend
to $\tet$ in the strong topology; in other words, $A$ is not
topologically bounded.
\end{proof}

As a consequence of the previous results, it is now possible to
characterize normability for characteristic \v Serstnev spaces
according to the following criterion.

\begin{theorem}\label{T:ch*S} A characteristic \v Serstnev space $(V,\nu,\tau)$ is
normable if, and only if, the null vector $\tet$ has a convex
$\D$--bounded neighbourhood.
\end{theorem}

The following (restrictive) sufficient condition is in
\cite{bP67}; we prove it here not only for the sake of
completeness, but also because Prochaska's thesis is not easily
accessible and, moreover, because the notation it adopts is
different from the one that has become usual after the publication
of \cite{book}.

\begin{theorem}[Prochaska]\label{T:pro}  A characteristic \v Serstnev space
$(V,\nu,\tau)$ with $\tau=\tau_M$ is locally convex.
\end{theorem}
\begin{proof} It suffices to consider the family of neighbourhoods of the origin $\tet$, $N_{\tet}(t)$, with $t>0$.
Let $t>0$, $p,q\in N_{\tet}(t)$ and $\al\in\lsp 0,1\rsp $. Then
\begin{equation*}
\begin{split}
\nu_{\al p+(1-\al)q}&(t)
\ge \tau_M\lp \nu_{\al p},\nu_{(1-\al)q}\rp (t)\\
 &=\sup_{\be\in\lsp 0,1\rsp }
 M\left(\nu_{\al p}(\be t),
 \nu_{(1-\al)q}\lp (1-\be)t\rp \right)\\
 &\ge M\left(\nu_{\al p}(\al t),
 \nu_{(1-\al)q}\lp (1-\al)t\rp\right) =
 M\left(\nu_p(t), \nu_q(t)\right)>1-t.
\end{split}
\end{equation*}
Thus $\al\,p+(1-\al)\,q$ belongs to $N_{\tet}(t)$ for every
$\al\in\lsp 0,1\rsp $.
\end{proof}

It is well known that every simple PN space $(V,\|\cdot\|,G,M)$
with $G\in\D^+$ satisfies the assumptions of Theorem \ref{T:pro}.
Moreover, these PN spaces are trivially normable, since their
strong topology coincides with the topology of their classical
norm. In general, it is expected that most of the PN spaces
considered in Theorem \ref{T:pro} will be normable, as shown by
the following corollary.

\begin{corollary}\label{C:cor2} Let $(V,\nu,\tau_M)$ be a characteristic \v Serstnev space. If $N_{\tet}(t)$ is
$\D$--bounded for some $t\in\lop 0,1\rop$, then $(V,\nu,\tau_M)$
is normable.
\end{corollary}

\subsection{Other cases}

Apart from the \v Serstnev spaces, we can also determine whether
an $\al$--simple space is normable, as the following result shows.

\begin{theorem} Let $G$ be a continuous and strictly increasing distribution function
in $\D^+$. Then, the $\al$--simple space $(V,\|\cdot\|,G;\al)$ is
normable.
\end{theorem}

\begin{proof} Let $N_{\tet}(t)$ be a neighbourhood of the origin $\tet$ with
$t\in\left] 0,1\right[$; then
\[
N_{\tet}(t)=\left\{p\in V : G\lp\frac t{\|
p\|^{\al}}\rp>1-t\right\}=\left\{p\in V : \| p\|<\lp\frac
t{G^{-1}(1-t)}\rp ^{1/\al}\right\}.
\]
Since $h(t)=\lp t/G^{-1}(1-t)\rp ^{1/\al}$ is a continuous
function such that $\lim_{t\to 0+}h(t)=0$ and $\lim_{t\to
1-}h(t)=\infty$, then it is clear that the strong topology for $V$
coincides with the topology of the norm $\| \cdot \|$ in $V$.
Therefore, $(V,\|\cdot\|,G;\al)$ is normable.
\end{proof}

It is natural to ask whether results similar to that of Theorem
\ref{T:ch*S} hold for general PN spaces. The conditions of Theorem
\ref{T:eq*S} need not be equivalent; for, there are PN spaces in
which a set $A$ may be topologically bounded without being
$\D$--bounded. On the other hand, even in those cases, sometimes
it is possible to establish directly whether a PN space that is
also a TV space is normable. To illustrate both facts, we next
introduce a new class of PN spaces whose interest excedes to serve
as an example at this point. Recall that only a few types of PN
spaces are known: finding a new type might be useful to deep into
the subject.

Before introducing the new class of  PN spaces we need the
following technical lemma.

\begin{lemma}\label{f} Let $f\colon\lsp
0,+\infty\rop\to\lsp 0,1\rsp$ be a right--continuous nonincreasing
function. Let we define $f^{[-1]}(1):=0$ and
$f^{[-1]}(y):=\sup\{x:f(x)>y\}$ for all $y\in [0,1[$
($f^{[-1]}(y)$ might be infinite). If $x_0\in \lsp 0,+\infty\rop$
and $y_0\in \lsp 0,1\rsp$, then the following facts are
equivalent: $(a)$ $f(x_0)>y_0; $ $(b)$ $x_0<f^{[-1]}(y_0)$.
\end{lemma}

\begin{proof} If $f(x_0)>y_0$ then
$f^{[-1]}(y_0)=\sup\{x:f(x)>y_0\}\ge x_0$. If we suppose that
$\sup\{x:f(x)>y_0\}= x_0$, then $f(x)\le y_0$ for every $x>x_0$.
Thus $f(x_0)=f(x_0+)\le y_0$, against the hypothesis; whence (a)
$\Rightarrow$ (b). The converse result is an immediate consequence
of the nonincreasingness of $f$.
\end{proof}

The following theorem introduces a new class of PN spaces---which
generalizes an example in \cite{bLG96}---, and also provides some
properties of the spaces in that class. As it has been said above,
such properties are interesting in order to our purposes in this
paper.

\begin{theorem}\label{T:f} Let $(V,\|\cdot\|)$ a normed space and let $T$ be a continuous $t$--norm. Let $f$ be a
function as in Lemma $\ref{f}$, and satisfying the following two
properties:
\begin{enumerate}
\item[(a)] $f(x)=1$ if and only if $x=0$;

\item[(b)] $f\left(\|p+q\|\right)\ge
T\left(f\left(\|p\|\right),f\left(\|q\|\right)\right)$ for every
$p,q\in V$.
\end{enumerate}
If $\nu\colon V\to\De^+$ is given by
\begin{equation}\label{E:nuf}
\nu_p(x)=\begin{cases} 0,&\quad x\le 0,\\ f(\|p\,\|), &\quad x\in\lop0,+\infty\rop,\\
             1, &\quad x=+\infty,\end{cases}
\end{equation}
for every $p\in V$, then $(V,\nu,\tau_T,\tau_{T^*})$ is a Menger
PN space satisfying the following properties:
\begin{enumerate}
\item[(F1)] $(V,\nu,\tau_T,\tau_{T^*})$ is a TV space;

\item[(F2)] $(V,\nu,\tau_T,\tau_{T^*})$ is normable;

\item[(F3)] If $p\in V$ and $t>0$, then the strong neighbourhood
$N_p(t)$ in $(V,\nu,\tau_T,\tau_{T^*})$ is not $\D$--bounded, but
$N_p(t)$ is topologically bounded whenever $N_p(t)\neq V$;

\item[(F4)] $(V,\nu,\tau_T,\tau_{T^*})$ is not a \v Serstev space;

\item[(F5)] $(V,\nu,\tau_T,\tau_{T^*})$ is not a characteristic PN
space.

\end{enumerate}
\end{theorem}

\begin{proof}
First, we prove that $(V,\nu,\tau_T,\tau_{T^*})$ satisfied the
four axioms to be a Menger PN space:

(N1) $\nu_p=\ep_0 \Leftrightarrow f(\|p\|)=1 \Leftrightarrow
\|p\|=0 \Leftrightarrow p=\tet$.\smallskip

(N2) Trivial.\smallskip

(N3) Given $p,q\in V$, then we have
$\nu_{p+q}\ge\tau_T\lp\nu_p,\nu_q\rp$ $\Leftrightarrow$
$\displaystyle\nu_{p+q}(x)\ge\tau_T\lp\nu_p,\nu_q\rp(x)=\sup_{s+t=x}
T\lp \nu_p(s),\nu_q(t)\rp$ for all $x\in\lop0,+\infty\rop$
$\Leftrightarrow$ $f(\|p+q\|)\ge T\lp f(\|p\|),f(\|q\|)\rp$, as
hyphothesized.\smallskip

(N4) Let $p\in V$ and let $\la\in \lsp 0,1\rsp$. Then,
$\nu_p\le\tau_{T^*}\lp\nu_{\la p},\nu_{(1-\la)p}\rp$
$\Leftrightarrow$ $\displaystyle\nu_p(x)\le\tau_{T^*}\lp\nu_{\la
p},\nu_{(1-\la)p}\rp(x)=\inf_{s+t=x} T^*\lp \nu_{\la
p}(s),\nu_{(1-\la)p}(t)\rp=\inf_{s+t=x} 1-T\lp 1-\nu_{\la
p}(s),1-\nu_{(1-\la)p}(t)\rp=1-\sup_{s+t=x} T\lp 1-\nu_{\la
p}(s),1-\nu_{(1-\la)p}(t)\rp$ for all $x\in\lop0,+\infty\rop$
$\Leftrightarrow$ $f(\|p\|)\le 1-\max
\{1-f(\la\|p\|),1-f((1-\la)\|p\|)\}=\min
\{f(\la\|p\|),f((1-\la)\|p\|)\}$. Therefore, for any $p\in V$,
$\nu_p\le\tau_{T^*}\lp\nu_{\la p},\nu_{(1-\la)p}\rp$ for all
$\la\in \lsp 0,1\rsp$ $\Leftrightarrow$ $f(\|p\|)\le f(\al\|p\|)$
for all $\al\in \lsp 0,1\rsp$. Hence,
$\nu_p\le\tau_{T^*}\lp\nu_{\la p},\nu_{(1-\la)p}\rp$ for all
$\la\in \lsp 0,1\rsp$ and for all $p\in V$ $\Leftrightarrow$ $f$
is nonincreasing.

Now we prove the four properties:

(F1) Let $p\in V$. We have to prove that the map from $\R$ into
$V$ defined by $\al\mapsto \al p$ is continuous at any $\al\in\R$.
Let $\gamma>0$ (we will suppose, without loss of generality, that
$\gamma\le 1$). We must prove that there exists a real number
$\de>0$ so that $d_S(\nu_{\al 'p-\al p},\ep_0)<\gamma$ whenever
$|\al '-\al|<\de$; or, equivalently, such that $d_S(\nu_{\be
p},\ep_0)<\gamma$ whenever $|\be|<\de$. Since
$d_S(\nu_q,\ep_0)=\inf\{h:\nu_q(h+)>1- h\}=1-f(\|q\|)$ (here
$\nu_q(h+)$ represents the limit $\lim_{u\to h^+}\nu_q(u)$), then
$d_S(\nu_{\be p},\ep_0)<\gamma \Leftrightarrow
1-f(|\be|\|p\|)<\gamma \Leftrightarrow f(|\be|\|p\|)>1-\gamma
\Leftrightarrow \textrm{ (Lemma 1) }
|\be|\|p\|<f^{[-1]}(1-\gamma)\Leftrightarrow
|\be|<f^{[-1]}(1-\gamma)/\|p\|=\de$.

(F2) Let $p\in V$. Let $t>0$ (we will suppose, without loss of
generality, that $t< 1-\lim_{x\to\infty}f(x)$). Then,
$N_p(t)=\{q\in V:d_S(\nu_{p-q},\ep_0)<t\} =\{q\in
V:1-f(\|p-q\|)<t\}=\{q\in V:f(\|p-q\|)>1-t\}=\textrm{ (Lemma 1) }
\{q\in V:\|p-q\|< f^{[-1]}(1-t)\}=B(p,f^{[-1]}(1-t))$, i.e., the
strong neighbourhood $N_p(t)$ is a ball in $(V,\|\cdot\|)$ with
center in $p$. Conversely, let $r>0$. If $t=1-f(r)$, then
$f^{[-1]}(1-t)<r$, whence $N_p(t)=B(p,f^{[-1]}(1-t))\subset
B(p,r)$. Therefore, the strong topology for
$(V,\nu,\tau_T,\tau_{T^*})$ coincides with the topology of the
norm in $(V,\|\cdot\|)$.

(F3) If $p\in V$ and $0<t< 1-\lim_{x\to\infty}f(x)$, then
$N_p(t)=B(p,f^{[-1]}(1-t))$ is a ball in $(V,\|\cdot\|)$, whence
$N_p(t)$ is topologically bounded. On the other hand, if
$0<x<\infty$ then $\Phi_{N_p(t)}(x)=\inf\{\nu_q(x):q\in
N_p(t)\}=\inf\{f(\|q\|):\|p-q\|<f^{[-1]}(1-t)\}=f\left(\|p\|+f^{[-1]}(1-t)\right)$.
Thus,
$\lim_{x\to\infty}R_{N_p(t)}(x)=f\left(\|p\|+f^{[-1]}(1-t)\right)<1$,
i.e.,  $N_p(t)$ is not $\D$--bounded.

(F4) It is immediate to check that $(V,\nu,\tau_T,\tau_{T^*})$ is
a \v Serstnev space if and only if the function $f$ is constant on
$]0,\infty[$. From hypothesis (a) this constant should be less
than 1, which contradicts the right--continuity of $f$ at $x=0$.
Thus, $(V,\nu,\tau_T,\tau_{T^*})$ is not a \v Serstnev space.

(F5) It is immediate that
$\nu(V\sm\{\tet\})\subseteq\De^+\sm\D^+$.
\end{proof}

Observe that the proof of (N4) is independent of the $t$--norm
$T$. Thus, we can also take $\tau^*=\tau_{M^*}=\tau_M$ (see
\cite{book}) in the PN--space of Theorem \ref{T:f}.

Now we consider some particular cases and provide some examples
which apply the preceding theorem.

\begin{example}\label{exPi}
Suppose that $T=\Pi$ in Theorem \ref{T:f}. In this case the
property (b) is read as $ f\lp \|p+q\|\rp \ge f\lp \|p\,\|\rp
\,f\lp \|q\|\rp$ for all $p,q\in V$. It is not difficult to prove
that, under the established hypotheses for $f$, the property (b)
is equivalent to the following one:
\begin{equation}\label{f:exPi}
f(x+y)\ge f(x)f(y)\,\, \textrm{for all }x,y\in [0,\infty[.
\end{equation}
It easy to check that instances of functions $f$ satisfying the
hypotheses of Theorem \ref{T:f} for this case are
\begin{gather*}
f_{\al,\be}(x):=1-\frac {\be}{\al}+\frac {\be}{x+\al},\quad 0\le
\be\le\al,\\
g_{\al,\be}(x):=1-\al+\al\exp\lp -x^{\be}\rp
,\quad0<\al\le1,\,\,\be>0.
\end{gather*}
\end{example}

\begin{example}\label{exW}
Suppose that $T=W$ in Theorem \ref{T:f}. In this case the property
(b) is read as $ f\lp \|p+q\|\rp \ge f\lp \|p\,\|\rp +f\lp
\|q\|\rp-1$ for all $p,q\in V$. Since $W$ is the minimum
continuous $t$--norm, all the functions $f$ satisfying the
hypotheses of Theorem \ref{T:f} with respect to any $t$--norm $T$
also satisfy such hypotheses with respect to $W$. It is not
difficult to prove that, under those hypotheses, the property (b)
is equivalent to the following one:
\[
1+f(x+y)\ge f(x)+f(y)\,\, \textrm{for all }x,y\in [0,\infty[.
\]
Instances of functions $f$ satisfying these hypotheses but not the
ones considered in Example \ref{exPi}---since they do not satisfy
(\ref{f:exPi})---are
\[
h_{\al,\be}(x):=\begin{cases} 1-\al x,&\quad 0\le x\le \be,\\
1-\al\be, &\quad x>\be,\end{cases}\quad 0<\be\le 1/\al.
\]
\end{example}

\section{Conclusion}

In what precedes we have been able to characterize those \v
Serstnev spaces that are normable TV spaces. Several questions
remain open: to give at least sufficient conditions under which a
general PN space is normable; more, to characterize (rather than
just having a sufficient condition) the class of PN spaces that
are also TV spaces, and, once this has been achieved, to study
normability in the class thus determined.

\end{document}